\documentclass[parskip=full]{llncs}
\usepackage{amssymb}
\usepackage{interval}
\usepackage{amsfonts}
\usepackage{hyperref}

\begin{document}

\newcommand{\overbar}[1]{\mkern
1.5mu\overline{\mkern-1.5mu#1\mkern-1.5mu}\mkern 1.5mu}

\newcommand{\realroot}[3]{Root(#1,[#2,#3])}

\newcommand{\isa}[1]{#1}

\title{Decidability of Univariate Real Algebra with Predicates for Rational and
Integer Powers}

\author{Grant Olney Passmore}
\institute{Aesthetic Integration, London and Clare Hall, University of
Cambridge\\
\email{grant.passmore@cl.cam.ac.uk}}

\maketitle

\begin{abstract}We prove decidability of univariate real algebra
  extended with predicates for rational and integer powers, i.e., ``$x^n \in
  \mathbb{Q}$'' and ``$x^n \in
  \mathbb{Z}$.'' Our decision procedure combines computation over real
  algebraic cells with the rational root theorem and witness
  construction via algebraic number density arguments.
  \end{abstract}

\section{Introduction}
From the perspective of decidability, the reals stand in stark
contrast to the rationals and integers.
While the elementary arithmetical theories of the integers and rationals are
undecidable, the corresponding theory of the reals is decidable and
admits quantifier elimination.
The immense utility real algebraic
reasoning finds within the mathematical sciences continues to motivate
significant progress towards practical automatic proof procedures for the
reals.

However, in mathematical practice, we are often faced with
problems involving a combination of nonlinear statements over the reals,
rationals
and integers.
Consider the existence and irrationality
of
$\sqrt{2}$, expressed in a language with variables implicitly ranging over
$\mathbb{R}$:
\[ \exists x (x \geq 0 \wedge x^2 = 2) \ \wedge \ \neg\exists x (x \in \mathbb{Q} \wedge x \geq 0
\wedge x^2 = 2) \]
Though easy to prove by hand this sentence has never to our knowledge
been placed within a broader decidable theory so that, e.g., the existence and
irrationality of solutions to any univariate real algebra problem can
be
decided
automatically.
This $\sqrt{2}$ example is relevant to the theorem proving
community
as its formalisation has been used as a benchmark for comparing
proof assistants~\cite{Wiedijk:2006:SPW:1121735}.
It would be useful if such proofs were fully automatic.

In this paper, we prove decidability of univariate real algebra extended
with
predicates
for rational and integer powers.
This guarantees we can always decide sentences like the above,
and many more besides.
For example, the following conjectures are decided by our method
in
a
fraction
of a second:
\[ \forall x (x^3 \in \mathbb{Z} \wedge x^5 \not\in \mathbb{Z} \ \Rightarrow \
x
\not\in \mathbb{Q}) \]
\[ \exists x (x^2 \in \mathbb{Q} \wedge x \not\in \mathbb{Q} \wedge x^5 + 1 >
20) \]
\[ \forall x (x^2 \not\in \mathbb{Q} \Rightarrow x \not\in \mathbb{Q}) \]
\[ \exists x (x \not\in \mathbb{Q} \wedge x^2 \in \mathbb{Z} \wedge 3x^4 + 2x +
1
>
5
\wedge
4x^3
+
1
<
2)
\]

\section{Preliminaries}\label{sec:prelim}

We assume a basic grounding in commutative algebra.
We do not however assume exposure to real algebraic geometry and give a
high-level treatment of the relevant
foundations.
%

The theory of real closed fields (RCF) is $Th(\langle \mathbb{R},
+,
-,
\times,
<,
0,
1
\rangle)$, the collection of all true sentences of the reals in the
elementary language of ordered rings.
RCF is complete, decidable and admits effective elimination of
quantifiers~\cite{Basu:2006:ARA:1197095}.

A real algebraic number is a real number that is a root of
 a (non-zero) univariate polynomial with integer coefficients.
 The real algebraic numbers, %
\[ \mathbb{R}_{alg} = \{x \in \mathbb{R} \ | \ \exists p \neq 0 \in
\mathbb{Z}[x] \ \textsf{s.t.} \ p(x) = 0\}, \]
form a computable subfield (a computable sub-RCF) of $\mathbb{R}$.
Indeed, $\mathbb{R}_{alg}$ embeds isomorphically into every RCF.
The field operations of $\mathbb{R}_{alg}$ are performed on computable
representations of field elements.
The \emph{minimal polynomial} of $\alpha \in \mathbb{R}_{alg}$ is the unique
monic $p \in \mathbb{Q}[x]$ of least degree s.t. $p(\alpha)=0$.
The degree of an algebraic number is the degree of its minimal polynomial.

An element $\alpha \in \mathbb{R}_{alg}$ can be represented by two pieces
of
data:
(i) a polynomial $p(x) \in \mathbb{Z}[x]$ s.t. $p(\alpha)=0$, and (ii) an
identifier specifying which root of $p(x)$ is denoted by $\alpha$.
A \emph{root-triple} representation is often used
where $\alpha$ is ``pinned down'' among the roots of $p(x)$ by an interval
with rational endpoints:
\[ \langle p(x) \in \mathbb{Z}[x], q_1, q_2 \in \mathbb{Q} \rangle \
\textsf{ s.t. } \
p(\alpha)=0 \ \wedge \ \#\{r \in [q_1, q_2] \ | \ p(r) = 0\}=1.\]

The process of \emph{root isolation} is a key component of
computing over $\mathbb{R}_{alg}$.
Given a polynomial $p \in \mathbb{Z}[x]$ with $k$ unique real roots,
root isolation computes a sequence of disjoint real
intervals with rational endpoints $I_1, \mathellipsis, I_k$ s.t.  each
$I_j$ contains precisely one real root of $p$.
Much work has been done on efficient root isolation.
Common approaches include those based on Sturm's Theorem and
Descartes' Rule of Signs \cite{collins1976polynomial,Mishra,uspensky:toe}.
Sturm's Theorem also plays a key role in computing the sign of a
polynomial evaluated at a real algebraic number.

Given representations of $\alpha, \beta \in \mathbb{R}_{alg}$,
there are two main approaches to performing the field operations, i.e., for
computing representations of ${\alpha}^{-1}$, $\alpha+\beta$, $\alpha\beta$,
etc.
Both approaches rely on root isolation.
The first approach uses bivariate resultants to compute representation
polynomials~\cite{Mishra}.
The second approach uses a recursive representation of real algebraic
numbers through an explicit treatment of field towers and does not require
computing resultants
\cite{Moura:2013aa,Rioboo:2003:TFR}.
Computing $\alpha^n$ (which plays a key role in our decision
procedure) can in general be done by repeated squaring, requiring on
the order of $\log n$ real algebraic number multiplications.  More
sophisticated methods for $\alpha^n$ are also available
\cite{Hirvensalo:PowersOfAlgNums}.

The Intermediate Value Theorem (IVT) holds over every RCF.
Armed with machinery for computing the sign of a polynomial $p(x) \in
\mathbb{Z}[x]$ at
a real algebraic
point $\alpha \in \mathbb{R}_{alg}$, the combination of IVT and root
isolation can be used as the basis of a decision method for univariate
real algebra.

\noindent Consider
\[ \varphi(x) = \left[\ \bigwedge_{i=1}^{k_1} \bigvee_{j=1}^{k_2}
\left(p_{ij}(x) \odot_{ij} 0\right) \right] \ \textsf{s.t.} \ p_{ij} \in
\mathbb{Z}[x], \ \odot_{ij} \in \{<, \leq, =, \geq, >\}. \]
We can decide the satisfiability of $\varphi$ over $\mathbb{R}$, i.e., whether
or not
\[ \langle \mathbb{R}, +, -, \times, <, 0, 1\rangle \ \models \ \exists
x(\varphi(x)) \]
in the following manner:
\begin{itemize}
 \item Let $P = \prod_{ij} p_{ij} \in \mathbb{Z}[x]$, the product of all
polynomials appearing in $\varphi$.
 \item Let $\alpha_1 < \mathellipsis < \alpha_k \in \mathbb{R}_{alg}$ be all
distinct real roots of $P$.
 \item Then, the roots $\alpha_i$ partition $\mathbb{R}$ into finitely many
connected components:
 \[ \mathbb{R} = \interval[open]{-\infty}{\alpha_1} \cup
                                 [\alpha_1] \cup
                                 \interval[open]{\alpha_1}{\alpha_2} \cup
                                 \dots \cup
                                 \interval[open]{\alpha_{k-1}}{\alpha_k} \cup
                                 [\alpha_k] \cup
                                 \interval[open]{\alpha_k}{+\infty}. \]
 \item By IVT, the sign of each polynomial $p_{ij}$ appearing in $\varphi$ is
invariant over any component of the partitioning.
 \item Thus, we can simply select one sample point from each component of the
partitioning and obtain a sequence of $2k+1$ real algebraic points $S = \{r_1,
\mathellipsis, r_{k+1}\} \subset \mathbb{R}_{alg}$ s.t.
\[ \langle \mathbb{R}, +, -, \times, <, 0, 1\rangle \ \models \ \exists
x(\varphi(x)) \ \ \iff \ \ \bigvee_{i=1}^{2k+1} \varphi(r_i).\]

Now $\exists x (\varphi(x))$ can be decided simply by evaluating $\varphi(x)$
at finitely many real algebraic points. The partitioning of $\mathbb{R}$
constructed above is called an \emph{algebraic decomposition} induced by $P$
(equivalently, by the polynomials $p_{ij}$).
\end{itemize}

\section{Decision Procedure}

Our decision procedure extends the IVT-based method for univariate real algebra
with means to handle predicates expressing the rationality and integrality of
powers of the variable of the formula, i.e., $(x^n \in \mathbb{Q})$ and $(x^n
\in \mathbb{Z})$.
As will be made clear (cf. Sec.~\ref{sec:discussion}), the restriction of these predicates to powers of the
variable is important: The method would fail if we allowed more general
polynomials $p(x) \in \mathbb{Z}[x]$ to appear in constraints of the form
$(p(x) \in \mathbb{Q})$.

Formally, we work over the univariate language of ordered rings $\mathcal{L}$
extended with infinitely many predicate symbols of one real variable:
\[(x \in \mathbb{Q}), (x^2 \in \mathbb{Q}), (x^3 \in \mathbb{Q}), \mathellipsis
\ \ \ \mbox{and} \ \ \ (x \in \mathbb{Z}), (x^2 \in \mathbb{Z}), (x^3 \in
\mathbb{Z}), \mathellipsis \ .
\]
We use $\mathcal{L}_{\mathbb{Q}\mathbb{Z}}$ to mean the resulting extended
language and $\mathcal{L}_{\mathbb{Q}}$ (resp. $\mathcal{L}_{\mathbb{Z}}$) to
mean $\mathcal{L}$ extended only with the rationality (resp. integrality)
predicates.

We present a method to decide the satisfiability of quantifier-free
$\mathcal{L}_{\mathbb{Q}\mathbb{Z}}$ formulas over $\mathbb{R}$.
It suffices to consider $\mathcal{L}_{\mathbb{Q}\mathbb{Z}}$ formulas of the
form
\[ \varphi(x) \wedge \Gamma(x) \]
where $\varphi \in \mathcal{L}$ is a formula of univariate real algebra and
\[\Gamma = \Gamma_{\mathbb{Q}} \wedge \Gamma_{\mathbb{Z}}\]
s.t.
\[ \Gamma_{\mathbb{Q}} =
   \left[\bigwedge_{i=1}^{k_1} (x^{w_{1}(i)} \in \mathbb{Q})
    \ \wedge \
    \bigwedge_{i=1}^{k_2} (x^{w_{2}(i)} \not\in \mathbb{Q})\right]\]
and
\[ \Gamma_{\mathbb{Z}} =
   \left[\bigwedge_{i=1}^{k_3} (x^{w_{3}(i)} \in \mathbb{Z})
    \ \wedge \
    \bigwedge_{i=1}^{k_4} (x^{w_{4}(i)} \not\in \mathbb{Z})\right].\]

Informed by the IVT-based method for univariate real algebra, we can reduce
this $\mathcal{L}_{\mathbb{Q}\mathbb{Z}}$ decision problem to an even more
restricted one. Crucial to this reduction is treating the connected components
of an algebraic decomposition as ``first class'' objects, rather than only
computing with single sample points selected from them.
We call such components \emph{r-cells}.

\begin{definition}[r-cell]
An r-cell is a connected component of $\mathbb{R}$ of one of the following four
forms (with $\alpha,\beta \in \mathbb{R}_{alg}$):
      (i) $[\alpha]$,
      (ii) $\interval[open]{-\infty}{\alpha}$ s.t. $\alpha \leq 0$,
      (iii) $\interval[open]{\alpha}{\beta}$ s.t. $0 \leq \alpha < \beta$ or
$\alpha < \beta \leq 0$,
      (iv) $\interval[open]{\alpha}{+\infty}$ s.t. $\alpha \geq 0$.
\end{definition}
Observe that the only r-cell containing zero is the singleton (type (i)) r-cell
$[0]$. Note that r-cells of type (i) are 0-dimensional subsets of $\mathbb{R}$
while r-cells of types (ii)-(iv) are 1-dimensional. We call these $0$-cells and
$1$-cells, resp.
An algebraic decomposition can always be transformed into an \emph{r-cell
decomposition} by splitting any 1-cell containing zero into three parts.

Given $\Phi(x) = \varphi(x) \wedge \Gamma(x)$, we must decide whether or not
$\mathbb{R}$ contains any point $x$ s.t. $\Phi(x)$ holds. To do so, we will
first compute an r-cell decomposition of $\mathbb{R}$ induced by the
polynomials of $\varphi$. Let $c_1, \mathellipsis, c_k$ be these r-cells. Then
by IVT, the truth of $\varphi$ is invariant within each $c_i$. Note,
however, that the truth of $\Gamma$ may vary over each $c_i$. Let
$C$
be
the
result of filtering out all r-cells $c_i$ that falsify $\varphi$:
\[ C = \{ c_i \ | \ \exists r \in c_i (\varphi(r)), \ 1 \leq i \leq k \}.\]
This can be done by evaluating $\varphi$ at a single sample point drawn from
each $c_i$.
If $C=\emptyset$, then $\Phi$ is clearly unsatisfiable over $\mathbb{R}$.
Otherwise, $C$ is a non-empty collection of r-cells over which $\varphi$ is
satisfied.
To decide $\Phi$, we need only to decide whether or not $\Gamma$ is satisfied
over any $c \in C$.

We present a method to do so.
We first develop a method to decide rationality constraints over an r-cell.
We then lift the method to handle general combinations of rationality and
integrality constraints.

\subsection{Deciding rationality constraints}

Given a system of rationality constraints $\Gamma_{\mathbb{Q}}$ and an r-cell
$c$, we need a method to decide whether or not $\Gamma_{\mathbb{Q}}$ is
satisfied over
$c$.
To accomplish this, we will extract a system of \emph{degree constraints} from
$\Gamma_{\mathbb{Q}}$ and give a method to decide if $c$ contains a real
algebraic number
satisfying them.

We must however take care of the following issue: If we prove there exists no
algebraic real in $c$ satisfying $\Gamma_{\mathbb{Q}}$, how do we know there
exists no
\emph{transcendental} real in $c$ satisfying $\Gamma_{\mathbb{Q}}$ as well?
That is, in the presence of rationality constraints, can we
still transfer results from $\mathbb{R}_{alg}$ to $\mathbb{R}$ as a whole?
We answer this question in the affirmative by proving a suitable transfer
principle
(cf.
Theorem~\ref{thm:transfer}).

It turns out we need essentially two methods for deciding
$\Gamma_{\mathbb{Q}}$ over $c$: One method for $0$-cells
and
another for $1$-cells.
We begin with the $1$-cell case.

\subsubsection{1-cells}

To construct our system of degree constraints, we shall
utilise
a
fundamental
property
relating
the
degree
of
a
``binomial root''
real algebraic number to the rationality of its powers.
We employ a result on the density of real algebraic numbers to show
that any consistent system of degree constraints gives rise to a real
algebraic solution in a 1-cell.
We then prove completeness of the method and a transfer principle enabling us
to
lift
results from $\mathbb{R}_{alg}$ to $\mathbb{R}$.

\begin{lemma}[Minimal binomials]\label{lem:irreducible-binomials}
Let $\alpha \in \mathbb{R}_{alg}$ s.t. $\alpha^n \in \mathbb{Q}$ for some $n
\in \mathbb{N}$.
Then,
the minimal polynomial for $\alpha$ over $\mathbb{Q}[x]$ is a binomial of the
form
$x^d
-
q$.
\begin{proof}
 Let $k \in \mathbb{N}$ be the least power s.t.
$\alpha^k
\in \mathbb{Q}$. We shall prove that $p(x) = x^k - \alpha^k \in \mathbb{Q}[x]$
is the minimal polynomial for $\alpha$.
Assume $p(x)$ is reducible over $\mathbb{Q}[x]$.
Observe that $p(x) = \prod_{i=1}^k (x - \alpha \zeta^i)$ where $\zeta$ is a
$kth$
root of unity.
As $p(x)$ is reducible, it must have a nontrivial factor $f(x) = \prod_{i=1}^m
(x - \alpha\zeta^{s_i}) \in \mathbb{Q}[x]$ with $m < k$ and $s_i \in
\mathbb{N}$.
But then $\left(\alpha^m \prod_{i=1}^m \zeta^{s_i}\right) \in \mathbb{Q}$, and
since $\alpha$ is real, we must have $\alpha^m \in \mathbb{Q}$.
But $m < k$. Contradiction. Thus, as $p(x) = x^k - \alpha^k$ is irreducible and
monic, it is the minimal polynomial for $\alpha$ over $\mathbb{Q}[x]$.
\qed
\end{proof}

\end{lemma}

\begin{lemma}[Binomial algebraic degree and divisibility]\label{lem:deg-div}
Let $\alpha \in \mathbb{R}_{alg}$ s.t. $\alpha$ is a root
of some $x^k - q \in \mathbb{Q}[x]$. Let $n \in
\mathbb{N}$.
Then,
\[ (\alpha^n \in \mathbb{Q}) \ \iff \ deg(\alpha) \mid n. \]
\begin{proof} Let $d = deg(\alpha)$.
$(\Leftarrow)$ By Lemma~\ref{lem:irreducible-binomials},
$\alpha^{d} \in \mathbb{Q}$. But, as $d \mid n$, we have
$\alpha^n = (\alpha^{d})^k$ for some $k \in \mathbb{N}$. Thus,
$\alpha^n
\in \mathbb{Q}$.
$(\Rightarrow)$ We use the method of infinite descent.
Consider $\alpha^n = q \in
\mathbb{Q}$.
Then, $x^n - q$ has $\alpha$ as a root, and thus $d \leq n$.
Assume $d \nmid n$.
It follows that $d<n$, $gcd(d,n)=1$, $q = \alpha^d \alpha^{n-d}$ and
$gcd(d,n-d)=1$.
As $\alpha^d \in \mathbb{Q}$, we have $\alpha^{n-d} = \frac{q}{\alpha^d} \in
\mathbb{Q}$.
Note $n-d < n$.
But then $\alpha^{n-d} \in \mathbb{Q}$ s.t. $d \nmid n-d$, and we can continue
this process ad
infinitum. Contradiction.
\qed
\end{proof}

\end{lemma}

Let $c \subset \mathbb{R}$ be a 1-cell and
$\Gamma_{\mathbb{Q}}$ a system of rationality constraints s.t.
\[\Gamma_{\mathbb{Q}}
=
\left[\bigwedge_{i=1}^{k_1} (x^{w_{1}(i)} \in \mathbb{Q})
  \ \wedge \
  \bigwedge_{i=1}^{k_2} (x^{w_{2}(i)} \not\in \mathbb{Q})\right].\]

\noindent To $\Gamma_{\mathbb{Q}}$, we associate a system of \emph{degree
constraints}
$\mathcal{D}(\Gamma_{\mathbb{Q}})$ as follows:

\[ \mathcal{D}(\Gamma_{\mathbb{Q}}) =
   \left[\bigwedge_{i=1}^{k_1} (d \mid w_1(i))
  \ \wedge \
  \bigwedge_{i=1}^{k_2} (d \nmid w_2(i))\right].\]

\noindent Note that each $w_j(i)$ is a concrete natural number.
Thus, $\mathcal{D}(\Gamma_{\mathbb{Q}})$ is a system of arithmetical
constraints with a
single
free
variable $d$. We shall prove that $\Gamma_{\mathbb{Q}}$ is satisfied over $c$
iff
$\mathcal{D}(\Gamma_{\mathbb{Q}})$ is consistent over $\mathbb{N}$, i.e., iff
\[\exists d \in \mathbb{N}
\textsf{
s.t.
}
\mathcal{D}(\Gamma_{\mathbb{Q}})(d).\]

We proceed in two steps.
First, we prove that $\Gamma_{\mathbb{Q}}$ is satisfied by a \emph{real
algebraic number} in
$c$ iff $\mathcal{D}(\Gamma_{\mathbb{Q}})$ is satisfied over $\mathbb{N}$.
Next, we show
that this result can be lifted to $\mathbb{R}$ as a whole, i.e., that
$\Gamma_{\mathbb{Q}}$
is satisfied over $c$ (by any real, be it algebraic or transcendental) iff
$\mathcal{D}(\Gamma_{\mathbb{Q}})$ is satisfied over $\mathbb{N}$.

These results elucidate a deep \emph{homogeneity} of $\mathbb{R}$.
Intuitively, $\mathbb{R}$ is so saturated with
real algebraic numbers that, given any open interval $I \subset \mathbb{R}$,
the only way $I$ can fail to contain an algebraic number satisfying
$\Gamma_{\mathbb{Q}}$ is
if the \emph{purely arithmetical facts} induced by $\Gamma_{\mathbb{Q}}$
(via Lemma~\ref{lem:deg-div})
are mutually inconsistent over $\mathbb{N}$. Moreover, from the perspective
of rationality constraints, transcendental elements cannot be distinguished
from
algebraic ones. To
prove
these results,
we
shall
need
to understand a bit about the density of real algebraic numbers of arbitrary
degree.

\begin{lemma}[Density of ratios of primes]\label{lemma:prime-density}
Given $a < b \in \mathbb{R}$, there exists $\frac{p}{q} \in
\interval[open]{a}{b}$ s.t. $|p|\not=|q|$ are both prime.
\begin{proof}
  A straightforward application of the Prime Number Theorem.
\end{proof}
\end{lemma}

\begin{lemma}[Density of real algebraic numbers of degree
n]\label{lemma:density}
  Let $a < b \in \mathbb{R}$ and $n \in \mathbb{N}$.  Then, $\exists
  \alpha \in \mathbb{R}_{alg}$
  s.t. $a < \alpha < b$ and $deg(\alpha) = n$ and $\alpha^n \in \mathbb{Q}$.
  \begin{proof}
    We construct an irreducible
     $p(x) = x^n - q \in \mathbb{Q}[x]$ s.t. $a <
    \sqrt[n]{q} < b$. Then, $\alpha = \sqrt[n]{q}$ will suffice.
    WLOG, assume $a>0$. Let $Q$ be a rational in
    $\interval[open]{a}{b}$. Let $f : \mathbb{R}^+ \to \mathbb{R}$ be
    the nth-root function, i.e., $f(r) = \sqrt[n]{r}$.
    Consider $Q^n \in \mathbb{Q}$. By continuity of $f$,
    $\exists\epsilon > 0$ s.t. $f(\interval[open]{Q^n -
      \epsilon}{Q^n + \epsilon}) \subset \interval[open]{a}{b}$.
    For each rational $q \in \interval[open]{Q^n - \epsilon}{Q^n + \epsilon}$,
    we thus have $a < f(q) < b$ with $f(q)$ algebraic, as $(f(q))^n - q = 0$.
    To prove the theorem, we must
    choose $q$ s.t. $deg(f(q)) = n$.
    It suffices to find $q \in \interval[open]{Q^n -
      \epsilon}{Q^n + \epsilon}$ s.t. $p(x) = x^n - q$ is irreducible
    over $\mathbb{Q}[x]$.
    By Lemma~\ref{lemma:prime-density}, we can choose $q = \frac{q_1}{q_2}
    \in \interval[open]{Q^n - \epsilon}{Q^n + \epsilon}$
    s.t. $q_1 \not= q_2$ are both prime.
    By Eisenstein's criterion, $q_2x^n - q_1$ is irreducible over
$\mathbb{Q}[x]$.
    Thus, $x^n - \frac{q_1}{q_2}$ is irreducible and
    $\alpha = \sqrt[n]{\frac{q_1}{q_2}}$ completes the proof.
    \qed
  \end{proof}
\end{lemma}

\noindent With Lemma~\ref{lemma:density} in hand, it is not hard to see that
$\Gamma_{\mathbb{Q}}$ is satisfied by a \emph{real algebraic number} in a
1-cell
$c$
iff
$\mathcal{D}(\Gamma_{\mathbb{Q}})$
is satisfied over $\mathbb{N}$.

\begin{theorem}[1-cell arithmetical reduction: algebraic
case]\label{thm:1dim-red-alg}
        Let $\Gamma_{\mathbb{Q}}$ be a system of rationality constraints and $c
\subseteq
\mathbb{R}$
a
1-cell. Then, $\Gamma_{\mathbb{Q}}$ is satisfiable over $c$ by a real
algebraic number iff
$\mathcal{D}(\Gamma_{\mathbb{Q}})$ is satisfiable over $\mathbb{N}$.
\begin{proof}
        $(\Rightarrow)$ Let $\alpha \in (c \cap \mathbb{R}_{alg})$ satisfy
$\Gamma_{\mathbb{Q}}$.
        Then, by Lemma~\ref{lem:deg-div}, $d = deg(\alpha)$ satisfies
$\mathcal{D}(\Gamma_{\mathbb{Q}})$.
        $(\Leftarrow)$ Let $d \in \mathbb{N}$ satisfy
$\mathcal{D}(\Gamma_{\mathbb{Q}})$.
        Then, by Lemma~\ref{lem:deg-div}, any algebraic $\alpha \in c$ s.t.
$deg(\alpha)=d$ will
satisfy $\Gamma_{\mathbb{Q}}$. But, by Lemma~\ref{lemma:density}, such an
$\alpha$ must
exist in $c$. \qed
\end{proof}

\end{theorem}

Thus, we have
reduced
the
satisfiability
of
$\Gamma_{\mathbb{Q}}$
by real algebraic numbers present in
a
1-cell $c$ to the satisfiability of $\mathcal{D}(\Gamma_{\mathbb{Q}})$ over
$\mathbb{N}$.
However, we must still attend to the possibility that $\Gamma_{\mathbb{Q}}$
could be
satisfied by a \emph{transcendental} element in $c$ without being satisfied by
an algebraic element in $c$. Let us now prove that this scenario is impossible.
In fact, we will prove this for both the 0 and 1-dimensional cases.

\begin{theorem}[Rationality constraints transfer principle]\label{thm:transfer}
Let $\Gamma_{\mathbb{Q}}$ be a system of rationality constraints and $c$ an
r-cell.
Then, it is impossible for $\Gamma_{\mathbb{Q}}$ to be satisfied by a
transcendental
real in $c$ without also being satisfied by an algebraic real in $c$.

\begin{proof}
Let
$\Gamma_{\mathbb{Q}} =
\left[\bigwedge_{i=1}^{k_1} (x^{w_{1}(i)} \in \mathbb{Q})
  \ \wedge \
  \bigwedge_{i=1}^{k_2} (x^{w_{2}(i)} \not\in \mathbb{Q})\right].$
  If $c$ is a 0-cell, then $c$ contains no transcendental elements,
  so the theorem holds.
  Consider $c$ a 1-cell.
  We examine the structure of $\Gamma_{\mathbb{Q}}$.
  If $k_1 > 0$, i.e., $\Gamma_{\mathbb{Q}}$ contains at least one positive
rationality
  constraint, then $\Gamma_{\mathbb{Q}}$ cannot be satisfied by any
transcendental element,
  and the theorem holds.
  Thus, we are left to consider $\Gamma_{\mathbb{Q}} = \bigwedge_{i=1}^{k_2}
  (x^{w_{2}(i)} \not\in \mathbb{Q})$ s.t. $\Gamma_{\mathbb{Q}}$ is satisfied by
a
  transcendental element in $c$.
  Let $m = max(w_2(1), \mathellipsis, w_2(k_2))$.
  Then, $\Gamma_{\mathbb{Q}}$ will be satisfied by any $\alpha \in \mathbb{R}_{alg}$ s.t.
  $deg(\alpha) > m$.
  But by Lemma~\ref{lemma:density}, $c$ must contain an algebraic $\alpha$
  s.t. $deg(\alpha) = m+1$.
  ~\qed
\end{proof}

\end{theorem}

In addition to giving us a complete method for deciding the satisfiability
of systems of rationality constraints over 1-cells, the
combination of
Theorem~\ref{thm:transfer} and the completeness of the theory of real closed
fields tells
us
something
of
a
fundamental
model-theoretic
nature:

\begin{corollary}[Transfer principle for $\mathcal{L}_{\mathbb{Q}}$]
        Given $\phi \in \mathcal{L}_{\mathbb{Q}}$,
        \[\langle \mathbb{R}, +, \times, <, (x^n \in \mathbb{Q})_{n \in
\mathbb{N}}
, 0, 1\rangle \models \phi
\\
\iff \\
        \langle \mathbb{R}_{alg}, +, \times, <, (x^n \in \mathbb{Q})_{n \in
\mathbb{N}}
, 0, 1\rangle \models \phi.\]
\end{corollary}

\noindent That is, extending the
language
$\mathcal{L}$
to
include
rationality constraints ($\mathcal{L}_{\mathbb{Q}}$) still guarantees a sound
transfer
of
results
from
$\mathbb{R}_{alg}$ to $\mathbb{R}$.

Finally, let us put the pieces together and prove our main
theorem for 1-cells.

\begin{theorem}[1-cell arithmetical reduction: general
case]\label{thm:transfer-general}
        Let $\Gamma_{\mathbb{Q}}$ be a system of rationality constraints and $c
\subseteq
\mathbb{R}$
a
1-cell. Then, $\Gamma_{\mathbb{Q}}$ is satisfiable over $c$ iff
$\mathcal{D}(\Gamma_{\mathbb{Q}})$ is satisfiable over $\mathbb{N}$.
\begin{proof}
Immediate by Theorem~\ref{thm:1dim-red-alg} and Theorem~\ref{thm:transfer}.
\qed
\end{proof}

\end{theorem}

Thus, to decide if $\Gamma_{\mathbb{Q}}$ is satisfied over a 1-cell $c$, we
need only check the consistency of $\mathcal{D}(\Gamma_{\mathbb{Q}})$ over
$\mathbb{N}$.
It is easy to derive an algorithm for doing so.
Consider $\mathcal{D}(\Gamma_{\mathbb{Q}})$ s.t. \[
\mathcal{D}(\Gamma_{\mathbb{Q}}) =
   \left[\bigwedge_{i=1}^{k_1} (d \mid w_1(i))
  \ \wedge \
  \bigwedge_{i=1}^{k_2} (d \nmid w_2(i))\right].\]

\noindent
If $k_1 = 0$, then $d = max(w_2(1), \mathellipsis, w_2(k_2))+1$
satisfies $\mathcal{D}(\Gamma_{\mathbb{Q}})$.
If $k_2 = 0$, then $d=1$ satisfies $\mathcal{D}(\Gamma_{\mathbb{Q}})$.
Finally, if $k_1>0$ and $k_2>0$, then $m = min(w_1(1), \mathellipsis,
w_1(k_1))$
gives us an upper bound on all $d$ satisfying
$\mathcal{D}(\Gamma_{\mathbb{Q}})$.
Thus, we need only search for such a $d$ from $1$ to $m$.
For efficiency, we can augment this bounded search by various cheap sufficient
conditions for recognising inconsistencies in
$\mathcal{D}(\Gamma_{\mathbb{Q}})$.

\subsubsection{0-cells}

When deciding rationality constraints over r-cells of the form $[\alpha]$,
we will need to decide, when given some $j \in \mathbb{N}$, whether or not
$\alpha^j \in \mathbb{Q}$.
Recall that a root-triple for $\alpha^j$ can be computed from a
root-triple for $\alpha$ (cf. Sec.~\ref{sec:prelim}).
A key component for deciding a system of rationality constraints over a
0-cell is
then
an
algorithm
for
deciding
whether
or
not
a
given real algebraic number $\beta = \alpha^j$ is rational.
Naively, one might try to solve this problem in the following way:

\begin{quote}
Given $\beta$ presented as a
root-triple
$\langle
p
\in \mathbb{Z}[x], l, u\rangle$, fully factor $p$ over $\mathbb{Q}[x]$.
Then, $\beta \in \mathbb{Q}$ iff the factorisation of $p$ contains a linear
factor of the form $(x - q)$ with $q \in \interval{l}{u}$.
\end{quote}

From the perspective of theorem proving, the problem with this approach
is that it is difficult in general to establish the ``completeness'' of a
factorisation.
While it is easy to verify that the product of a collection of factors
equals the original polynomial, it can be very challenging (without
direct appeal to the functional correctness of an implemented
factorisation algorithm) to prove that a given polynomial is
irreducible, i.e., that it cannot be factored any further.
Indeed, deep results in algebraic number theory are used even to
classify the irreducible factors of binomials \cite{Hollmann1986}.
Moreover, univariate factorisation can be computationally expensive,
especially when one is only after rational roots.

We would like the steps in our proofs to be as clear and obvious as possible,
and to minimise the burden of formalising our procedure as a tactic in a proof
assistant.
Thus, we shall go a different route.
To decide whether or not a given $\alpha$ is rational, we apply a
simple but powerful result from high school mathematics: %

\begin{theorem}[Rational roots]\label{thm:rat-roots}
  Let $p(x) = \sum_{i=0}^n a_nx^n \in \mathbb{Z}[x] \setminus \{0\}$.
  If $\frac{a}{b} \in \mathbb{Q}$ s.t. $p(q) = 0$ and $gcd(a,b)=1$,
  then $a \mid a_0$ and $b \mid a_n$.
  \begin{proof}
    A straightforward application of Gauss's lemma.
  \end{proof}
\end{theorem}

Given Theorem~\ref{thm:rat-roots}, we can decide the rationality of $\alpha$
simply by enumerating potential rational roots $q_1,
\mathellipsis, q_k$ and checking by evaluation whether any $q_i$ satisfies
$\left(l \leq
q_i
\leq r \ \wedge \ p(q_i)=0\right)$.
Then, to decide whether $\alpha$ satisfies a given system of
rationality constraints,
e.g.,
$\Gamma_{\mathbb{Q}} = \left[(x^2
\in \mathbb{Q}) \wedge (x \not\in \mathbb{Q})\right]$, we first compute
a root-triple representation for $\alpha^2$ and
then
test $\alpha$ and $\alpha^2$ for rationality as described.
This process clearly always terminates.
To make this more efficient when faced with many potential rational
roots,
we
can
combine
(i)
dividing
our
polynomial
$p$
by
$(x-q)$ whenever $q$ is realised to be a rational root, and (ii)
various
cheap
irreducibility criteria over $\mathbb{Q}[x]$
for recognising when a polynomial
has no linear factors over $\mathbb{Q}[x]$ and thus has no rational roots.

\subsection{Deciding integrality constraints}

\subsubsection{Integrality constraints over an unbounded 1-cell}

\newcommand{\clgamma}[0]{\overbar{\Gamma}}

WLOG let $c = \interval[open]{\alpha}{+\infty}$ with $\alpha \geq 0$.
Consider $\Gamma = \Gamma_{\mathbb{Q}} \wedge \Gamma_{\mathbb{Z}}$ with
\[ \Gamma_{\mathbb{Z}} =
   \left[\bigwedge_{i=1}^{k_3} (x^{w_{3}(i)} \in \mathbb{Z})
    \ \wedge \
    \bigwedge_{i=1}^{k_4} (x^{w_{4}(i)} \not\in \mathbb{Z})\right].\]
    We use the notation $\phi : \Gamma$ to mean that the constraint $\phi$ is
present as a conjunct in $\Gamma$.
It is convenient to also view $\Gamma$ as a set.
Let $\clgamma$ denote the closure of $\Gamma$
under
the
following
saturation rules:

\begin{enumerate}
 \item $(x^n \not\in \mathbb{Q}) : \clgamma \ \rightarrow \ (x^n
\not\in \mathbb{Z}) : \clgamma$
 \item $(x^n \in \mathbb{Z}) : \clgamma \ \rightarrow \ (x^n \in
   \mathbb{Q}) : \clgamma$
 \item $(x^n \in \mathbb{Z}) : \clgamma \wedge (x^m \not\in
\mathbb{Z}) : \clgamma
\
\rightarrow \
(x \not\in \mathbb{Q}) : \clgamma$
 \item $(x^n \in \mathbb{Z}) : \clgamma \wedge (x^m \in \mathbb{Q})
:
\clgamma \ \rightarrow \
(x^m \in \mathbb{Z}) : \clgamma$
 \item $(x^n \in \mathbb{Z}) : \clgamma \wedge (x^m \not\in \mathbb{Z}) :
\clgamma \
\rightarrow \ (x^m \not\in \mathbb{Q}) : \clgamma$
\end{enumerate}

\noindent
This saturation process is clearly finite.
The
soundness
of
rules
1 and 2
is
obvious.
The soundness of rules 3-5 is easily verified by the following lemmata.

\begin{lemma}[Soundness: rule 3]\label{lem:rule-3}
        $(x^n \in \mathbb{Z}) \wedge (x^m \not\in \mathbb{Z})
   \rightarrow (x \not\in \mathbb{Q})$
   \begin{proof}
        Since $x^m \not\in \mathbb{Z}$, we know $x \not\in \mathbb{Z}$.
        Suppose $x \in \mathbb{Q}$. Then $x = \frac{a}{b}$ s.t. $gcd(a,b)=1$.
        Thus, $a^n = x^n b^n$.
        Thus, $b \mid a$.
        Recall $gcd(a,b)=1$. So, $b=1$.
        But then $x = a \in \mathbb{Z}$.
        Contradiction.
        \qed
   \end{proof}

\end{lemma}

\begin{lemma}[Soundness: rule 4]\label{lem:rule-4}
        $(x^n \in \mathbb{Z}) \wedge (x^m \in \mathbb{Q})
   \rightarrow (x^m \in \mathbb{Z})$
   \begin{proof}
     Let $d = deg(x)$. By Lemma~\ref{lem:deg-div}, $d \mid n$ and $d \mid m$.
     If $d = n$, then $x^m = (x^n)^k$ for some $k \in \mathbb{N}$ and thus $
x^m \in \mathbb{Z}$.
     Otherwise, $d < n$. Let $x^{d} = \frac{a}{b} \in \mathbb{Q}$ s.t.
$gcd(a,b)=1$.
     Thus, $x^n = (x^d)^k = \frac{a^k}{b^k} \in \mathbb{Z}$ for some
$k \in \mathbb{N}$.
     But then $b = 1$, and thus $x^d \in \mathbb{Z}$.
     So, as $d \mid m$, $x^m \in \mathbb{Z}$ as well.
     \qed
   \end{proof}

\end{lemma}

\begin{lemma}[Soundness: rule 5]\label{lem:rule-5}
        $(x^n \in \mathbb{Z}) \wedge (x^m \not\in \mathbb{Z})
   \rightarrow (x^m \not\in \mathbb{Q})$
   \begin{proof}
     Assume $(x^n \in \mathbb{Z})$ and $(x^m \not\in \mathbb{Z})$ but $(x^m \in
\mathbb{Q})$.
     But then $(x^m \in \mathbb{Z})$ by rule 4. Contradiction.
     \qed
   \end{proof}

\end{lemma}

Let us now prove
that
these
rules\footnote{In fact, the completeness proof shows that rule 3 is
logically unnecessary. Nevertheless, we find its inclusion in the saturation process
useful
in
practice.}
are \emph{complete}
for
deciding the satisfiability of systems of rationality and integrality
constraints over unbounded 1-cells. Let $\clgamma_{\mathbb{Q}}$
(resp. $\clgamma_{\mathbb{Z}}$)
 denote the collection of rationality (resp. integrality) constraints present
in
$\clgamma$.
Intuitively,
we
shall
exploit
the
following
observation: The construction of $\clgamma$
projects all information pertaining to the \emph{consistency} of the
combined
rationality
and
integrality
constraints of $\Gamma$
onto
$\clgamma_{\mathbb{Q}}$. Then, if $\clgamma_{\mathbb{Q}}$
is
consistent, i.e., $\exists d \in \mathbb{N}$ satisfying
$\mathcal{D}(\clgamma_{\mathbb{Q}})$, this will impose a strict
correspondence between $\clgamma_{\mathbb{Q}}$ and
$\clgamma_{\mathbb{Z}}$. From this correspondence and a least
$d$ witnessing $\mathcal{D}(\clgamma_{\mathbb{Q}})$,
we
can
construct an algebraic real satisfying $\Gamma$.

\begin{lemma}[$\clgamma_{\mathbb{Q}}$-$\clgamma_{\mathbb{Z}}$
correspondence]\label{lem:q-z-correspondence}
If $\Gamma_{\mathbb{Z}}$ contains at least one positive
integrality constraint, then
 \[ \forall m \in \mathbb{N} \left[(x^m \in \mathbb{Q}) : \clgamma
\
\iff \ (x^m \in \mathbb{Z}) : \clgamma\right] \]
and
 \[ \forall m \in \mathbb{N} \left[(x^m \not\in \mathbb{Q}) : \clgamma
\
\iff \ (x^m \not\in \mathbb{Z}) : \clgamma\right]. \]

\begin{proof}
 Let us call the first conjunct A and the second B.
 $(A\Rightarrow)$ As
$\Gamma_{\mathbb{Z}}$ contains at least one positive integrality constraint,
rule 4 guarantees $(x^m \in \mathbb{Z}) : \clgamma$.
 $(A\Leftarrow)$ Immediate by rule 2.
$(B\Rightarrow)$ Immediate by rule 1.
$(B\Leftarrow)$ As $\Gamma_{\mathbb{Z}}$ contains at least one positive
integrality constraint, rule 5 guarantees $(x^m \not\in \mathbb{Q}) :
\clgamma$.
\qed
\end{proof}
\end{lemma}

\begin{theorem}[Completeness of $\Gamma$-saturation
method]\label{thm:gamma-sat-ub-complete}
Let $\Gamma = \Gamma_{\mathbb{Q}} \wedge \Gamma_{\mathbb{Z}}$ be a system of
rationality and integrality constraints, and
$c
\subseteq
\mathbb{R}$ an unbounded 1-cell.
Then, $\mathcal{D}(\clgamma_{\mathbb{Q}})$ is consistent over
$\mathbb{N}$ iff $\Gamma$ is consistent over $c$.
\begin{proof}
$(\Leftarrow)$
Immediate by Theorem~\ref{thm:transfer-general} and the soundness of our
saturation rules.
$(\Rightarrow)$
We proceed by cases.

[Case 1: $\Gamma$ contains no positive rationality constraint]:
Then, by Lemma~\ref{lem:q-z-correspondence} and the consistency of
$\mathcal{D}(\clgamma_{\mathbb{Q}})$, $\Gamma_{\mathbb{Z}}$ must contain no
positive integrality constraints. But then it is consistent with $\Gamma$ that
every
power of $x$ listed in $\Gamma$ be irrational. Let $k \in \mathbb{N}$ be
the largest power s.t. $x^k$ appears in a constraint in $\Gamma$. Then, by
Lemma~\ref{lem:deg-div}, any $\alpha \in c$ s.t. $deg(\alpha)>k$ will satisfy
$\Gamma$. By Lemma~\ref{lemma:density}, we can always find such an $\alpha$ in
$c$, e.g., we can select $\alpha \in c$ s.t. $deg(\alpha) = k+1$.

[Case 2: $\Gamma$ contains a positive rationality constraint but no positive
integrality constraints]: By the
consistency of $\mathcal{D}(\clgamma_\mathbb{Q})$, it is consistent with
$\Gamma$
for
every
power
of
$x$ listed in $\Gamma$ to be non-integral. Let $d \in \mathbb{N}$ be the least
natural number satisfying $\mathcal{D}(\Gamma_{\mathbb{Q}})$. Then, we can
satisfy
$\Gamma$ with an $\alpha$ s.t. $deg(\alpha) = d$ with $\alpha^{dk} \in
\left(\mathbb{Q} \setminus
\mathbb{Z}\right)$ for each $x^{dk}$ appearing in a constraint in $\Gamma$.
By Lemma~\ref{lemma:density}, we know such an $\alpha$ is present in $c$ of the
form $\alpha = \sqrt[d]{\frac{p}{q}}$ for primes $p \not= q$.

[Case 3: $\Gamma$ contains both positive rationality and integrality
constraints]
By Lemma~\ref{lem:q-z-correspondence}, the rows of $\clgamma_{\mathbb{Q}}$ and
$\clgamma_{\mathbb{Z}}$ are in perfect correspondence.
Let $d \in \mathbb{N}$ be the least natural number satisfying
$\mathcal{D}(\clgamma_{\mathbb{Q}})$.
Since $\clgamma_{\mathbb{Q}}$ is consistent, we can satisfy $\Gamma$ by
finding an $\alpha \in c$ s.t. $\alpha^{dk} \in \mathbb{Z}$ for every $x^{dk}$
appearing in a constraint in $\Gamma$.
Recall $c$ is unbounded towards $+\infty$.
Thus, $c$ contains infinitely many primes $p$ s.t. $\sqrt[d]{p} \in c$.
Let $p \in c$ be such a prime.
Then, $x^d - p \in \mathbb{Q}[x]$ is irreducible by Eisenstein's criterion.
Thus, $\sqrt[d]{p} \in c$ and satisfies $\Gamma$.
\qed
\end{proof}

\end{theorem}

\subsubsection{Integrality constraints over a bounded 1-cell}

Let us now consider the satisfiability of $\Gamma = \Gamma_{\mathbb{Q}}
\wedge \Gamma_{\mathbb{Z}}$ over a bounded 1-cell $c \subset
\mathbb{R}$.
Given the results of the last section, it is easy to see that
if $\mathcal{D}(\clgamma_{\mathbb{Q}})$ is unsatisfiable over $\mathbb{N}$,
then $\Gamma$ is unsatisfiable over $c$.
However, as $\Gamma$ is bounded on both sides, it is possible for
$\mathcal{D}(\clgamma_{\mathbb{Q}})$ to be satisfiable over
$\mathbb{N}$ while $\Gamma$ is unsatisfiable over $c$.
That is, provided $\mathcal{D}(\clgamma_{\mathbb{Q}})$ is consistent
over $\mathbb{N}$, we must find a way to determine if $c$ actually
contains some $\alpha$ s.t. $\Gamma(\alpha)$ holds.
Afterall, even with $\clgamma_{\mathbb{Q}}$ satisfied over $c$, it is
possible that $c$ itself is not ``wide enough'' to satisfy the
integrality constraints $\clgamma_{\mathbb{Z}}$.
WLOG, let $c = \interval[open]{\alpha}{\beta}$ s.t. $0 \leq
\alpha < \beta \in \mathbb{R}_{alg}$.
Let $\mathcal{D}(\clgamma_{\mathbb{Q}})$ be satisfied by $d \in
\mathbb{N}$.
If $\Gamma$ contains no positive integrality
constraints, then we can reason as we did in the proof of
Theorem~\ref{thm:gamma-sat-ub-complete} to show $\Gamma$ is satisfied
over $c$.
The difficulty arises when a positive constraint $(x^k \in
\mathbb{Z})$ appears in $\Gamma_{\mathbb{Z}}$.
We can solve this case as follows.

\begin{theorem}[Satisfiability over a bounded
1-cell]\label{thm:sat-bounded}
Let $\Gamma_{\mathbb{Z}}$ contain at least one positive integrality constraint.
Let $\mathcal{D}(\clgamma_{\mathbb{Q}})$ be satisfiable over $\mathbb{N}$ with $d \in \mathbb{N}$
the least witness.
Let $c = \interval[open]{\alpha}{\beta}$ s.t. $0 \leq \alpha < \beta \in
\mathbb{R}_{alg}$.
  Then, $\Gamma$ is satisfiable over $c$
iff $\exists z \in \left(\interval[open]{\alpha^d}{\beta^d} \cap
\mathbb{Z}\right)$ s.t. $x^d - z \in \mathbb{Z}[x]$ is irreducible over
$\mathbb{Q}[x]$.
  \begin{proof}
    $(\Rightarrow)$
    Assume $\Gamma$ is satisfied by $\alpha \in c$.
Then, by soundness of $\clgamma$ saturation, $\clgamma$ is satisfied by
$\alpha$ as well.
    By Lemma~\ref{lem:q-z-correspondence}, $(x^d \in \mathbb{Z}) : \clgamma$.
    Moreover, $d$ is the least natural number with this property.
As $0 \leq \alpha < \beta$, $\{r^d \ | \ r \in c\} =
\interval[open]{\alpha^d}{\beta^d}$.
Thus, as $\Gamma$ is satisfied by $\alpha \in c$, there must exist an integer
$z \in \interval[open]{\alpha^d}{\beta^d}$ s.t. $deg(\sqrt[d]{z}) = d$.
    But then by uniqueness of minimal polynomials, $x^d - z$ is irreducible
over
$\mathbb{Q}[x]$.
    $(\Leftarrow)$
Assume $z \in \left(\interval[open]{\alpha^d}{\beta^d} \cap \mathbb{Z}\right)$
s.t. $x^d - z$ is irreducible over $\mathbb{Q}[x]$.
Let $\gamma = \sqrt[d]{z}$ and note that $\gamma \in
\interval[open]{\alpha}{\beta}.$
By Lemma~\ref{lem:deg-div}, $deg(\gamma) = d$.
Thus, $\clgamma_{\mathbb{Q}}$ is satisfied by $\gamma$.
As $\gamma^d \in
\mathbb{Z}$, it follows by Lemma~\ref{lem:q-z-correspondence} that
$\Gamma$
is
satisfied by $\gamma$ as well.
    \qed
  \end{proof}

\end{theorem}

\noindent By Eisenstein's criterion, we obtain a useful corollary.

\begin{corollary}
Let $\mathcal{D}(\clgamma_{\mathbb{Q}})$ be satisfiable with $d \in \mathbb{N}$
the least natural number witness.
Let $c = \interval[open]{\alpha}{\beta}$ s.t. $0 \leq \alpha < \beta \in
\mathbb{R}_{alg}$.
  Then, $\Gamma$ is satisfiable over $c$ if
  $\exists p \in \interval[open]{\alpha^d}{\beta^d}$ s.t. $p$ is prime.

\end{corollary}

These results give us a simple algorithm to decide satisfiability of
$\Gamma$ over $c$:
If $\mathcal{D}(\clgamma_{\mathbb{Q}})$ is unsatisfiable over $\mathbb{N}$,
then
$\Gamma$ is unsatisfiable. Otherwise, let $d \in \mathbb{N}$ be the minimal
solution
to
$\mathcal{D}(\clgamma_{\mathbb{Q}})$.
Gather all integers $\{z_1, \mathellipsis, z_k\}$ in $I =
\interval[open]{\alpha^d}{\beta^d}$.
If any $z_i$ is prime, $\Gamma$ is satisfied over $c$.
Otherwise, for each $z_i$, form the real
algebraic number
$\sqrt[d]{z_i}$ and check by evaluation if it satisfies
$\Gamma$.
By Theorem~\ref{thm:sat-bounded},
$\Gamma$ is satisfiable over $c$ iff one of the $\sqrt[d]{z_i} \in c$
satisfies this process.

\subsubsection{Integrality constraints over a 0-cell}

Finally, we consider the case of $\Gamma = \Gamma_{\mathbb{Q}} \wedge
\Gamma_{\mathbb{Z}}$ over a 0-cell $[\alpha]$.
Clearly, $\Gamma$ is satisfied over $c$ iff $\Gamma$ is satisfied at $\alpha$.
By the soundness of $\Gamma$-saturation, if
$\mathcal{D}(\clgamma_{\mathbb{Q}})$
is
unsatisfiable over $\mathbb{N}$, then $\Gamma$ is unsatisfiable over $c$.
Thus, we first form $\clgamma$ and check satisfiability of
$\mathcal{D}(\clgamma_{\mathbb{Q}})$ over $\mathbb{N}$.
Provided it is satisfiable, we then check $\Gamma(x \mapsto \alpha)$ by
evaluation.

\section{Examples}

We have implemented\footnote{The implementation of our procedure, including
computations over r-cells,
$\Gamma$-saturation and the proof output routines can be found
in
the
{\tt
RCF/}
modules
in
the MetiTarski source code at \url{http://metitarski.googlecode.com/}.} our
decision
method
in
a
special
version
of
the
MetiTarski
theorem prover \cite{paulson2012metitarski}.
We do not use any of the proof search mechanisms of MetiTarski, but rather
its parsing and first-order formula data structures.

In the examples that follow, all output (including the prose and \LaTeX\space
formatting)
has
been
generated
automatically by our implementation of the method.

\subsection{Example 1}
Let us decide $\exists x (\varphi(x) \wedge \Gamma(x))$, where
 \[ \varphi = \isa{(x ^ 2 - 2 = 0)} \ \textnormal{ and } \
   \Gamma = (x \in \mathbb{Q}).\]

\noindent We first compute $\clgamma$, the closure of $\Gamma$ under the
saturation rules:
  \[\clgamma = (x \in \mathbb{Q}).\]
  Observe $\mathcal{D}(\clgamma_{\mathbb{Q}})$ is satisfied (minimally) by
$d=1$.

We next compute an r-cell decomposition of $\mathbb{R}$ induced by $\varphi$,
yielding:
\begin{enumerate}
 \item $\interval[open]{-\infty}{\realroot{x^2 -2}{-2}{-1/3}}$,
 \item $[\realroot{x^2 -2}{-2}{-1/3}]$,
 \item $\interval[open]{\realroot{x^2 -2}{-2}{-1/3}}{0}$,
 \item $[0]$,
 \item $\interval[open]{0}{\realroot{x^2 -2}{1/3}{2}}$,
 \item $[\realroot{x^2 -2}{1/3}{2}]$,
 \item $\interval[open]{\realroot{x^2 -2}{1/3}{2}}{+\infty}$.
\end{enumerate}

\noindent By IVT, $\varphi$ has constant truth value over each such r-cell.
Only two r-cells in the decomposition satisfy $\varphi$:

  $[\realroot{x^2 -2}{-2}{-1/3}]$,
  $[\realroot{x^2 -2}{1/3}{2}]$.

\noindent Let us now see if any of these r-cells satisfy $\Gamma$.
\begin{enumerate}
 \item We check if $[\realroot{x^2 -2}{-2}{-1/3}]$ satisfies $\Gamma$.
 \begin{enumerate}
  \item Evaluating $(\alpha \in \mathbb{Q})$ for $\alpha=\realroot{x^2
-2}{-2}{-1/3}$.
  We shall determine the numerical type of $\alpha$.
  Let $p(x) = x^2 -2$.
     By RRT and the root interval, we reduce the set of possible rational
values for $\alpha$ to
      $\{-1, -2\}$.
     But none of these are roots of $p(x)$.
  Thus, $\alpha \in (\mathbb{R}\setminus\mathbb{Q})$.
 \end{enumerate}
 So, the r-cell does not satisfy $\Gamma$.
 \item We check if $[\realroot{x^2 -2}{1/3}{2}]$ satisfies $\Gamma$.
 \begin{enumerate}
  \item Evaluating $(\alpha \in \mathbb{Q})$ for $\alpha=\realroot{x^2
-2}{1/3}{2}$.
  We shall determine the numerical type of $\alpha$.
  Let $p(x) = x^2 -2$.
     By RRT and the root interval, we reduce the set of possible rational
values for $\alpha$ to
      $\{1, 2\}$.
     But none of these are roots of $p(x)$.
  Thus, $\alpha \in (\mathbb{R}\setminus\mathbb{Q})$.
 \end{enumerate}
 So, the r-cell does not satisfy $\Gamma$.
\end{enumerate}
Thus, as all r-cells have been ruled out, the conjecture is false. \qed

\subsection{Example 2}
Let us decide $\exists x (\varphi(x) \wedge \Gamma(x))$, where
 \[ \varphi = \isa{True} \ \textnormal{ and } \
   \Gamma = (x^3 \in \mathbb{Z}) \wedge (x^5 \not\in \mathbb{Z}) \wedge (x \in
\mathbb{Q}).\]

\noindent We first compute $\clgamma$, the closure of $\Gamma$ under the
saturation rules:
  \[\clgamma = (x \not\in \mathbb{Z}) \wedge (x \in \mathbb{Z}) \wedge (x
\not\in \mathbb{Q}) \wedge (x \in \mathbb{Q}) \wedge (x^3 \in \mathbb{Q})
\wedge (x^3 \in \mathbb{Z}) \wedge (x^5 \not\in \mathbb{Q}) \wedge (x^5 \not\in
\mathbb{Z}).\]
But, $\clgamma$ is obviously inconsistent.
Thus, the conjecture is false. \qed

\subsection{Example 3}
Let us decide $\exists x (\varphi(x) \wedge \Gamma(x))$, where
 \[ \varphi = \isa{((x ^ 3 - 7 > 3) \wedge (x ^ 2 + x + 1 < 50))} \ \textnormal{ and } \
   \Gamma = (x^2 \not\in \mathbb{Q}) \wedge (x^3 \in \mathbb{Z}).\]

\noindent We first compute $\clgamma$, the closure of $\Gamma$ under the saturation rules:
  \[\clgamma = (x^2 \not\in \mathbb{Q}) \wedge (x^2 \not\in \mathbb{Z}) \wedge (x^3 \in \mathbb{Z}) \wedge (x^3 \in \mathbb{Q}).\]
  Observe $\mathcal{D}(\clgamma_{\mathbb{Q}})$ is satisfied (minimally) by $d=3$.

We next compute an r-cell decomposition of $\mathbb{R}$ induced by $\varphi$, yielding:
\begin{enumerate}
 \item $\interval[open]{-\infty}{\realroot{x^2 + x -49}{-8}{-1/50}}$,
 \item $[\realroot{x^2 + x -49}{-8}{-1/50}]$,
 \item $\interval[open]{\realroot{x^2 + x -49}{-8}{-1/50}}{0}$,
 \item $[0]$,
 \item $\interval[open]{0}{\realroot{x^3 -10}{57/44}{5/2}}$,
 \item $[\realroot{x^3 -10}{57/44}{5/2}]$,
 \item $\interval[open]{\realroot{x^3 -10}{57/44}{5/2}}{\realroot{x^2 + x -49}{401/100}{8}}$,
 \item $[\realroot{x^2 + x -49}{401/100}{8}]$,
 \item $\interval[open]{\realroot{x^2 + x -49}{401/100}{8}}{+\infty}$.
\end{enumerate}

\noindent By IVT, $\varphi$ has constant truth value over each such r-cell. Only one r-cell in the decomposition  satisfies $\varphi$:

  $\interval[open]{\realroot{x^3 -10}{57/44}{5/2}}{\realroot{x^2 + x -49}{401/100}{8}}$.

\noindent Let us now see if any of these r-cells satisfy $\Gamma$.
\begin{enumerate}
 \item We check if $\interval[open]{\realroot{x^3 -10}{57/44}{5/2}}{\realroot{x^2 + x -49}{401/100}{8}}$ satisfies $\Gamma$.
   Call the boundaries of this r-cell $L$ and $U$.
 As $\Gamma$ contains a positive integrality constraint and $d=3$, any satisfying witness in this r-cell must be of the form $\sqrt[3]{z}$ for z an integer in $\interval[open]{L^3}{U^3}$.
 The set of integers in question is $Z=\{ z \in \mathbb{Z} \ | \ 11 \leq z \leq 276\}$, containing 266 members.
We shall examine $\sqrt[3]{z}$ for each $z \in Z$ in turn.
 \begin{enumerate}
  \item Evaluating $(\alpha^2 \not\in \mathbb{Q})$ for $\alpha=\realroot{x^3 -11}{1/12}{11}$.
   Observe $\alpha^2 = \realroot{x^3 -121}{1/144}{121}$.
  We shall determine the numerical type of $\alpha^2$.
  Let $p(x) = x^3 -121$.
     By RRT and the root interval, we reduce the set of possible rational values for $\alpha^2$ to
      $\{1, 11, 121\}$.
     But none of these are roots of $p(x)$.
  Thus, $\alpha^2 \in (\mathbb{R}\setminus\mathbb{Q})$.
  \item Evaluating $(\alpha^3 \in \mathbb{Z})$ for $\alpha=\realroot{x^3 -11}{1/12}{11}$.
   Observe $\alpha^3 = \realroot{x^3 -1331}{1/1728}{1331}$.
  We shall determine the numerical type of $\alpha^3$.
  Let $p(x) = x^3 -1331$.
     By RRT and the root interval, we reduce the set of possible rational values for $\alpha^3$ to
      $\{1, 11, 121, 1331\}$.
  Thus, we see $\alpha^3 =11 \in \mathbb{Z}$.
 \end{enumerate}
Witness found: $\realroot{x^3 -11}{1/12}{11}$.
 So, the r-cell does satisfy $\Gamma$.
\end{enumerate}
Thus, the conjecture is true. \qed

\section{Discussion and Related Work}\label{sec:discussion}

Let us describe some related results that help put our work into context.
\begin{itemize}
 \item The existence of rational or integer solutions to univariate polynomial
equations over $\mathbb{Q}[x]$ has long been known to be decidable. The best
known algorithms are based on univariate
factorisation via lattice reduction \cite{vanHoeij2002167}.
 \item Due to Weispfenning, the theory of linear, multivariate mixed
real-integer
arithmetic
is known to be decidable and admit quantifier elimination~\cite{Weispfenning:1999:MRL:309831.309888}.
 \item Due to van den Dries, the theory of real closed fields extended with a predicate for powers of two is known to be decidable~\cite{vandendries:rcf2}. Avigad and Yin have given a syntactic decidability proof for this theory, establishing a non-elementary upper bound for eliminating a block of quantifiers~\cite{Avigad200748}.
 \item Due to Davis, Putnam, Robinson and Matiyasevich, the $\exists^3$
nonlinear,
equational
theories
of
arithmetic over $\mathbb{N}$ and $\mathbb{Z}$ are known to be undecidable
(``Hilbert's Tenth Problem'' and reductions of its negative
solution)~\cite{matiyasevich:Book}.
 \item The decidability of the $\exists^2$ nonlinear,
equational theories of arithmetic over $\mathbb{N}$ and $\mathbb{Z}$ is open.
 \item Due to Poonen, the $\forall^2\exists^7$ theory of nonlinear arithmetic
over
$\mathbb{Q}$ is known to be undecidable \cite{Poonen}. This is an improvement
of
Julia
Robinson's original undecidability proof of $Th(\mathbb{Q})$ via a
$\forall^2\exists^7\forall^6$
definition
of
$\mathbb{Z}$
over
$\mathbb{Q}$~\cite{robinson:phd}.
 \item Due to Koenigsmann, the $\forall^{418}$ and $\forall^1\exists^{1109}$ theories of nonlinear arithmetic
over
$\mathbb{Q}$ are known to be
undecidable, via explicit definitions of $\mathbb{Z}$ over
$\mathbb{Q}$~\cite{koenigsmann,koenigsmann:personal}.
 \item The decidability of the $\exists^k$
equational nonlinear theory of arithmetic over
$\mathbb{Q}$ is open for $k>1$ (``Hilbert's Tenth Problem over $\mathbb{Q}$'').
\end{itemize}
Our present result --- the decidability of the nonlinear, univariate
theory of the reals extended with predicates for rational and integer
powers --- fills a gap somewhere between the positive result on
linear, multivariate mixed real-integer arithmetic, and the
negative result for Hilbert's Tenth Problem in three variables.

Next, we would like to turn our decision method
into a verified proof procedure within a proof assistant.
The deepest result needed is the Prime Number Theorem (PNT).
As Avigad et al have formalised a proof of PNT
within
Isabelle/HOL
\cite{Avigad:2007:FVP:1297658.1297660}, we are hopeful that a verified version
of our procedure can be built in Isabelle/HOL~\cite{paulson1994isabelle} in the near future.
To this end, it is useful to observe that PNT is not needed by
the restriction of our method to deciding the
rationality
of
real algebraic numbers like $\sqrt{2}$ and ${\sqrt{3}+\sqrt{5}}$.
Thus, a simpler tactic could be constructed for this
fragment.

Finally, we hope to extend the method to allow constraints of the form $(p(x) \in \mathbb{Q})$ for more general polynomials $p(x) \in \mathbb{Z}[x]$.
The key difficulty lies with Lemma~\ref{lem:deg-div}.
This crucial property relating the degree of an algebraic number to the rationality of its powers applies to ``binomial root'' algebraic numbers, but not to algebraic numbers in general.
For example, consider $\alpha = \sqrt{2}+\sqrt[4]{2}$. Then, the minimal polynomial of $\alpha$ over $\mathbb{Q}[x]$ is $x^4 - 4x^2 - 8x + 2$, but $\alpha^4 \not\in \mathbb{Q}$.
Thus, in the presence of richer forms of rationality and integrality constraints, our degree constraint reasoning is no longer sufficient.
We expect to need more powerful tools from algebraic number theory to extend the method in this way.

\section{Conclusion}

We have established decidability of univariate real algebra extended with
predicates for rational and integer powers.
Our decision procedure combines computations over real algebraic cells with the
rational root theorem and results on the density of real
algebraic numbers.
We have implemented the method, instrumenting it to produce readable proofs.
In the future, we hope to extend our result to richer systems of rationality
and integrality constraints, and to construct a verified version of the
procedure within a proof assistant.

\paragraph*{Acknowledgements.}
{We thank Jeremy Avigad, Wenda Li, Larry Paulson, Andr\'as Salamon and the anonymous referees for their helpful comments.}

\bibliographystyle{splncs03}

\bibliography{real_closed_field}

\end{document}